\documentclass{amsart}

\newtheorem{theorem}{Theorem}[section]
\newtheorem{lemma}[theorem]{Lemma}
\newtheorem{corollary}[theorem]{Corollary}

\theoremstyle{definition}

\theoremstyle{remark}

\usepackage{amsfonts, amssymb, latexsym}

\numberwithin{equation}{section}

\newcommand{\C}{\mathbb{C}}

\begin{document}

\title
[$K$-Theory of Crepant Resolutions of Complex  $SU(2)$-Orbifolds]
{
$K$-Theory of Crepant Resolutions of Complex Orbifolds with $SU(2)$
Singularities
}

\author{Christopher Seaton
}

\address{Department of Mathematics and Computer Science,
Rhodes College, 2000 N. Parkway,
Memphis, TN 38112}

\email{seatonc@rhodes.edu}

\subjclass{Primary 19L47}

\date{May 19, 2005}

\keywords{Orbifold, orbifold crepant resolution, orbifold $K$-theory}

\begin{abstract}
We show that if $Q$ is a closed, reduced, complex orbifold of dimension $n$ such that
every local group acts as a subgroup
of $SU(2) < SU(n)$, then the $K$-theory of the unique crepant resolution of $Q$
is isomorphic to the orbifold $K$-theory of $Q$.

\end{abstract}

\maketitle


\section{Introduction}

Let $Q$ be a reduced, compact, complex orbifold of dimension $n$; i.e. a compact
Hausdorff space locally modeled on $\C^n/G$ where $G$ is a finite group which acts
effectively on $\C^n$ with a fixed-point set of codimension at least 2
(for details of the definition and further background, see \cite{ruangwt}).  Then a
crepant resolution of $Q$ is given by a pair $( Y, \pi )$ where $Y$ is a smooth complex
manifold of dimension $n$ and $\pi : Y \rightarrow Q$ is a surjective map which is
biholomorphic away from the singular set of $Q$, such that $\pi^\ast K_Q = K_Y$ where
$K_Q$ and $K_Y$ denote the canonical line bundles of $Q$ and $Y$, respectively (see \cite{joyce}
for details).  In \cite{ruansgt}, it is conjectured that if $\pi : Y \rightarrow Q$ is a
crepant resolution of a Gorenstein orbifold $Q$ (i.e. an orbifold such that all groups
act as subgroups of $SU(n)$), then the orbifold $K$-theory of $Q$
is isomorphic to the ordinary $K$-theory of $Y$.  For the case of a global quotient
of $\C^n$, this has been verified for $n=2$ in \cite{nakajima} and, for Abelian
groups and a specific choice of crepant resolution for $n=3$ in \cite{itonakajima}.
Here, we apply the `local' results in the case $n=2$ to the case of a general orbifold with
such singularities.

The $K$-theory of an orbifold can be defined in several different ways.
First, it can be defined in the usual way in terms of equivalence classes of
orbifold vector bundles (see \cite{ademruan}).  As well, it is well-known
that a reduced orbifold $Q$ can be expressed as the quotient $P/G$ where $P$ is a
smooth manifold and $G$ is a compact Lie group \cite{kawasaki2}.
In the case of a real orbifold, $P$ can be taken to be the orthonormal
frame bundle of $Q$ with respect to a Riemannian metric and $G = O(n)$.  Similarly,
in the complex case, $P$ can be taken to be the unitary frame bundle and $G = U(n)$.
Hence, the orbifold $K$-theory of $Q$ is defined as the $G$-equivariant $K$-theory
$K_G(P)$.  See \cite{ademruan} or \cite{marcolli} for more details.

In Section \ref{decomp}, we describe the structure of the singular set $\Sigma$ of $Q$ in the
case in question and state the main result.  In section \ref{thproof}, we interpret
this decomposition in terms of ideals of the $C^\ast$-algebra of $Q$ and prove the result.

The author would like to thank Carla Farsi, Yongbin Ruan, and Siye Wu
for helpful suggestions and discussions leading to this result.
This work was supported in part by the Thron Fellowship of the University
of Colorado Mathematics Department.


\section{The Decomposition of $\Sigma$ and Statement of the Result}
\label{decomp}

Let $Q$ be a closed, reduced, complex orbifold with $\mbox{dim}_{\C}\, Q = n$, and fix a hermitian metric
on $TQ$ throughout.  Then each point $p \in Q$ is contained in a neighborhood modeled by $\C^n/G_p$ where
$p$ corresponds to the origin in $\C^n$ and $G_p < U(n)$.
Suppose that each of the local groups $G_p$ act as a subgroup of $SU(2) < SU(n)$, and then
each point $p$ is locally modeled by $\C^n / G_p \cong \C^{n - 2} \times (\C^2 / G_p)$.  Suppose further that $Q$ admits a crepant
resolution $\pi : Y \rightarrow Q$ so that $Y$ is a closed complex $n$-manifold.  By Proposition 9.1.4 of
\cite{joyce}, $(Y, \pi)$ is a {\bf local product resolution}, which in this context means the following
(see 9.1.2 of \cite{joyce} for the general definition):

Fix $p \in Q$, and then there is a neighborhood $U_p \ni p$ modeled by $\C^n / G_p$.  By hypothesis,
$U_p \cong V \times W/G_p$
where $V \times \{ 0 \} \cong \C^{n - 2}$ is the fixed point set of $G_p$, $W \cong \C^2$ is the orthogonal
complement of $V$ in $\C^n$ (for some choice of $G_p$-invariant metric on $\C^n$), and we identify
$G_p < SU(n)$ with its restriction $G_p < SU(2)$.  Then for a resolution
$( Y_p, \pi_p )$ of $W/G_p$, we let $\phi : V \times W/G_p \rightarrow \C^n / G_p$, $T$ be the ball of radius
$R > 0$ about the origin in $\C^n / G_p$, and $U := (\mbox{id} \times \pi_p )^{-1} (T) \subset V \times Y_p$.
There is a local isomorphism $\psi : (V \times Y_p) \backslash U \rightarrow Y$ such that the following
diagram commutes
\[
\begin{array}{ccc}
    (V \times Y_p) \backslash U
        &   \stackrel{\psi}{\longrightarrow}    &   Y   \\\\
        \downarrow \mbox{id} \times \pi_p       &&
        \downarrow \pi                          \\\\
    (V \times W/G_p) \backslash T
        &   \stackrel{\phi}{\longrightarrow}    &   \C^n/G_p .
\end{array}
\]
Hence, each of the singular points in a neighborhood of $p$ is resolved by $V \times Y_p$.
Moreover, as $(Y, \pi)$ is a crepant resolution of $Q$, $(Y_p, \pi_p )$ is a crepant resolution
of $\C^2/G_p$ (\cite{joyce}, Proposition 9.1.5), and hence is the unique crepant resolution
of $\C^2/G_p$.  It is clear that a crepant resolution of $Q$ can be formed by patching together
local products of the unique crepant resolutions of $\C^2 / G_p$, but we now see that this
is the only crepant resolution of $Q$.  Moreover, if $S$ denotes a connected component
of the singular set $\Sigma$ of $Q$, then a neighborhood of $S$ can be covered by a finite
number of charts as above, so that the isotropy subgroups of any $p, q \in S$ are conjugate
in $SU(2)$.  Moreover, each such chart $\C^n / G^p \cong V \times W/G_p$ restricts to a
complex manifold chart of dimension $n-2$ for $S$.

We summarize this discussion in the following.

\begin{lemma}
\label{decomplemma}
Let $Q$ be a closed, reduced, complex orbifold of complex dimension $n$, and suppose each of the
local groups $G_p$ acts on $Q$ as a subgroup of $SU(2)$.  Then there is a unique crepant resolution
$(Y, \pi )$ of $Q$.  The singular set $\Sigma$ of $Q$ is given by
\[
    \Sigma = \bigsqcup\limits_{i=1}^k S_i
\]
for some $k$ finite, where each $S_i$ is a connected, closed, complex $(n-2)$-manifold and the (conjugacy
class of the) isotropy subgroup $G_p < SU(2)$ of $p$ is constant on $S_i$.  Moreover, if $N_i$ is a sufficiently
small tubular neighborhood of $S_i$ in $Q$, then $N_i \cong S_i \times \C^2/G_p$ and
$\pi^{-1} (N_i) \cong S^i \times Y_i$ where $Y_i$ is the unique crepant resolution of $\C^2 / G_p$.

\end{lemma}

Such a decomposition may be possible for orbifolds with $SU(3)$ singularities; in this case,
components of the singular set have $(n-2)$- and $(n-3)$-dimensional components.  The latter
are closed manifolds, but the former may be open.
However, the techniques in this paper do not easily extend to this case.  For finite
subgroups of $SU(3)$, crepant resolutions are not unique.  While a local isomorphism
has been constructed for abelian subgroups of $SU(3)$ (see \cite{itonakajima}), this
is for a specific choice of resolution.

Using the decomposition given in this lemma, we will show the following:

\begin{theorem}
\label{mnthrm}
Let $Q$ be a closed, reduced, complex orbifold of complex dimension $n$, and suppose each of the
local groups $G_p$ acts on $Q$ as a subgroup of $SU(2)$.  Let $(Y, \pi )$ denote the unique crepant
resolution of $Q$, and then
\[
    K_{orb}^\ast (Q) \cong K^\ast (Y)
\]
as additive groups.

\end{theorem}

For any $n$-dimensional orbifold that admits a crepant resolution, the local groups can be chosen
to be subgroups of $SU(n)$ (see \cite{joyce}).  Therefore, we have as an immediate corollary:

\begin{corollary}
\label{dim2cor}
Let $Q$ be a 2-dimensional complex orbifold which admits a crepant resolution $(Y, \pi )$.
Then
\[
    K_{orb}^\ast (Q) \cong K^\ast (Y)
\]
as additive groups.

\end{corollary}


\section{Proof of Theorem \ref{mnthrm}}
\label{thproof}

In order to prove Theorem \ref{mnthrm}, we will show that $K_\ast (A) \cong K_\ast (B)$
where $A$ is the $C^\ast$-algebra of $Q$ and $B$ the $C^\ast$-algebra of $Y$.
So fix an orbifold $Q$ that satisfies the hypotheses of Theorem \ref{mnthrm}, and let
$k$, $S_i$, $N_i$, etc. be as given in Lemma \ref{decomplemma}.
We assume that the $N_i$ are chosen small enough so that $N_i \cap N_j = \emptyset$
for $i \neq j$.

For each $i$, let $N_i^\prime$ be a smaller tubular neighborhood
of $S_i$ so that $S_i \subset N_i^\prime \subset \overline{N_i^\prime} \subset N_i$,
and let $N_0 := Q \backslash \bigcup\limits_{i=1}^k \overline{N_i^\prime}$.  Then
$\{ N_i \}_{i = 0}^k$ is an open cover of $Q$ such that $N_0$ contains no singular
points.  Note that the restriction $\pi_{| \pi^{-1} (N_0)}$ is a biholomorphism
onto $N_0$.

Let $P$ denote the unitary frame bundle of $Q$, and then $Q = P / U(n)$.
Let $A: = C^\ast(Q)$ denote the $C^\ast$-algebra
$C(P) \rtimes_\alpha U(n)$ of $Q$ where $\alpha$ is the action of $U(n)$ on $C(P)$
induced by the usual action on $P$, and let $A^0$ denote the dense subalgebra
$L^1(U(n), C(P), \alpha)$ of $C(P) \rtimes_\alpha U(n)$.  Let $I_1^0$
denote the ideal in $A^0$ consisting of functions $\phi$ such that
$\phi (g)$ vanishes on $P_{|S_1}$ for each $g \in U(n)$
(i.e $I_1^0 = L^1 (U(n), C_0(P \backslash P_{|S_1}), \alpha )$;
as usual, $P_{|S_1}$ denotes the restriction of $P$ to $S_1$),
and let $I_1$ be the closure of $I_1^0$ in $A$.  Similarly, for each $j$ with
$1 < j \leq k$, set
$I_j^0 := L^1 \left( U(n), C_0 \left( P \backslash \bigcup\limits_{i=1}^j P_{|S_i} \right) , \alpha \right)$
to be the ideal of functions $\phi$ in $A^0$ such that for each $g \in U(n)$, $\phi(g)$ vanishes on the fibers over
$S_1, S_2, \cdots , S_j$, and $I_j$ the closure of $I_j^0$ in $A$.  Then we have
the ideals
\[
    I_k \subset I_{k-1} \subset \cdots \subset I_1 \subset I_0 := A .
\]
Note that, for each $j$ with $1 \leq j < k$,
$I_j/I_{j+1} \cong C(P_{|S_{j+1}}) \rtimes_\alpha U(n)$, and
$I_k \cong C_0(P_{|N_0}) \rtimes_{\alpha} U(n)$.

Similarly, let $B := C(Y)$ denote the algebra of continuous functions on $Y$,
and let $J_j$ denote the ideal of functions which vanish on
$\pi^{-1} \left( \bigcup\limits_{i=1}^j S_i \right)$.  Then we have
\[
    J_k \subset J_{k-1} \subset \cdots \subset J_1 \subset J_0 :=B ,
\]
with $J_j / J_{j+1} \cong C(\pi^{-1}(S_{j+1}))$ and $J_k \cong C_0(\pi^{-1}(N_0) )$.

Recall that $\pi$ restricts to a biholomorphism
\[
    \pi_{|\pi^{-1}(N_0)} : \pi^{-1}(N_0)
    \stackrel{\cong}{\longrightarrow}
    N_0 .
\]
Hence, as the action of $U(n)$ is free on $P_{|N_0}$,
\[
\begin{array}{rcl}
    K_\ast (I_k)
    &=&     K_\ast(C_0 (P_{|N_0}) \rtimes_\alpha U(n) ) \\\\
    &\cong&     K_{U(n)}^\ast (P_{|N_0})                \\\\
    &&\mbox{naturally, by the Green-Julg Theorem (\cite{black} Theorems 20.2.7 and 11.7.1)},
                                    \\\\
    &\cong&     K^\ast (P_{|N_0} /U(n) )                \\\\
    &&\mbox{as the $U(n)$ action is free on $N_0$},         \\\\
    &=&     K^\ast (N_0)                    \\\\
    &=&     K^\ast (\pi^{-1}(N_0) )             \\\\
    &=&     K_\ast (J_k) .
\end{array}
\]
Therefore, there is a natural isomorphism
\begin{eqnarray}
    K_\ast (I_k)
    &\cong&     K_\ast (J_k) .
\end{eqnarray}
Hence, Theorem \ref{mnthrm} holds for orbifolds such that $k = 0$; i.e. manifolds.
The next lemma gives an inductive step which, along with the previous result,
yields the theorem.

\begin{lemma}
\label{inductivestep}
Suppose
\[
    K_\ast (I_j ) \cong K_\ast (J_j )
\]
naturally for some $j$ with $1 \leq j \leq k$.  Then
\[
    K_\ast (I_{j-1} ) \cong K_\ast (J_{j-1} ) .
\]

\end{lemma}

\begin{proof}

Note that $I_j$ is an ideal in $I_{j-1}$, with $I_{j-1}/I_j = C(P_{|S_j}) \rtimes_\alpha U(n)$.
Similarly, $J_j$ is an ideal in $J_{j-1}$ with $J_{j-1}/J_j = C(\pi^{-1}(S_j))$.
We have the standard exact sequences
\[
\begin{array}{ccccc}
    K_0 (I_j)   &   \rightarrow &   K_0 (I_{j-1})   &   \rightarrow &   K_0 (I_{j-1}/I_j)   \\\\
    \partial\uparrow
            &           &           &           &   \downarrow \partial \\\\
    K_1(I_{j-1}/I_j)&   \leftarrow  &   K_1 (I_{j-1})   &   \leftarrow  &   K_1 (I_j)
\end{array}
\]
and
\[
\begin{array}{ccccc}
    K_0 (J_j)   &   \rightarrow &   K_0 (J_{j-1})   &   \rightarrow &   K_0 (J_{j-1}/J_j)   \\\\
    \partial\uparrow
            &           &           &           &   \downarrow \partial \\\\
    K_1(J_{j-1}/J_j)&   \leftarrow  &   K_1 (J_{j-1})   &   \leftarrow  &   K_1 (J_j) .
\end{array}
\]
So if we show that $K_\ast (I_{j-1}/I_j) \cong K_\ast (J_{j-1}/J_j)$ naturally, by the Five lemma, we are done.

Note that $I_{j-1}/I_j$ is the $C^\ast$-algebra of the quotient orbifold $P_{|S_j}/U(n)$,
which is given by the smooth manifold $S_j$ with the trivial action of $G_j$ (here, $G_j$ denotes a
choice from the conjugacy class of isotropy groups $G_p$ for $p \in S_j$).  Hence,
$I_{j-1}/I_j \cong C(S_j) \otimes C^\ast(G_j)$.
Similarly, we have
\[
\begin{array}{rcl}
    J_{j-1}/J_j
    &=& C(\pi^{-1}(S_j))                \\\\
    &=& C(S_j \times Y_j)               \\\\
    &=& C(S_j) \otimes C(Y_j),
\end{array}
\]
where $Y_j$ is the preimage of the origin in the unique crepant resolution of $\C^2/G_j$.  However,
$K_0(C^\ast (G_j)) = R(G)$ (\cite{black} Proposition 11.1.1 and Corollary 11.1.2)
which is naturally isomorphic to $K^0(Y_j)$ by \cite{nakajima} (Section 4.3; see also
\cite{itonakajima}), and $K^0(Y_j) \cong K_0 (C(Y_j))$, so that $K_0 (C^\ast(G_j))$
and $K_0 (C(Y_j))$ are isomorphic.  With this, by the K\"{u}nneth Theorem for
tensor products (\cite{black} Theorem 23.1.3),
\[
\begin{array}{cccccc}
    0   &\rightarrow&   K_0(C(S_j))\otimes K_0(C^\ast(G_j)) &\rightarrow& K_0(C(S_j) \otimes C^\ast(G_j))
        &\rightarrow                        \\\\
        &   &   \downarrow              &   &   \downarrow      \\\\
    0   &\rightarrow&   K_0(C(S_j))\otimes K_0(C(Y_j))  &\rightarrow&   K_0(C(S_j)\otimes C (Y_j))
        &\rightarrow
\end{array}
\]
\bigskip
\bigskip
\[
\begin{array}{cccc}
    \rightarrow&    \mbox{Tor}(K_0(C(S_j)), K_0(C^\ast(G_j)))           &\rightarrow&   0           \\\\
        &   \downarrow                                      \\\\
    \rightarrow&    \mbox{Tor}(K_0(C(S_j)), K_0(C(Y_j)))    &\rightarrow&   0
\end{array}
\]
and the Five lemma, we have a natural isomorphism
\[
    K_0 (C(S_j) \otimes C^\ast (G_j)) \cong K_0 (C(S_j) \otimes C (Y_j)) .
\]
So
\[
    K_0 (I_{j-1}/I_j) \cong K_0 (J_{j-1}/J_j).
\]

For the $K_1$ groups, we note that by \cite{black}, Corollary 11.1.2, $K_1(C^\ast (G_j)) = 0$.
As well, $K_1 (C(Y_j)) \cong K^1(Y_j)$, and it is known that $Y_j$ is diffeomorphic to a finite
collection of 2-spheres which intersect at most transversally at one point (see \cite{joyce}).
Therefore, $K^1(Y_j) = 0$.  Here, the hypothesis that all groups act as subgroups of $SU(2)$
is crucial.  For subgroups of $SU(3)$, the topology of the resolution is
not understood
sufficiently to compute the $K_1$ groups.

With this, we again apply the K\"unneth theorem and Five lemma
\[
\scriptstyle
\begin{array}{ccccccc}
    0   &\rightarrow&   K_1(C(S_j))\otimes K_0(C^\ast(G_j)) \oplus
                K_0(C(S_j))\otimes K_1(C^\ast(G_j)) &\rightarrow&   K_1(C(S_j) \otimes C^\ast(G_j))
        &&
                                                        \\\\
        &&      \downarrow          &&      \downarrow          \\\\
    0   &\rightarrow&   K_1(C(S_j))\otimes K_0(C(Y_j)) \oplus K_0(C(S_j))\otimes K_1(C(Y_j))
                &\rightarrow&   K_1(C(S_j) \otimes C(Y_j))  &&
\end{array}
\]
\bigskip
\bigskip
\[
\begin{array}{ccccc}
        \cdots &\rightarrow&    \mbox{Tor}(K_1(C(S_j)), K_0(C^\ast(G_j))) \oplus
                \mbox{Tor}(K_0(C(S_j)), K_1(C^\ast(G_j)))   &\rightarrow&   0       \\\\
        &&      \downarrow                                  \\\\
        \cdots & \rightarrow&   \mbox{Tor}(K_1(C(S_j)), K_0(C(Y_j))) \oplus
                \mbox{Tor}(K_0(C(S_j)), K_1(C(Y_j)))
        \rightarrow&    0
\end{array}
\]
Therefore, we have a natural isomorphism
\[
    K_1 (C(S_j) \otimes C^\ast (G_j)) \cong K_1 (C(S_j) \otimes C(Y_j)) ,
\]
and
\[
    K_1 (I_{j-1}/I_j) \cong K_1 (J_{j-1}/J_j).
\]

\end{proof}

Now, as $K_\ast (I_k) \cong K_\ast (J_k)$, repeated application of Lemma \ref{inductivestep}
yields that $K_\ast (A) \cong K_\ast (B)$, and hence we have proven Theorem \ref{mnthrm}.


\bibliographystyle{amsplain}

\end{document}